\documentstyle[twoside,12pt]{article}
\newtheorem{prop}{Proposition}[section]
\newtheorem{dfn}[prop]{Definition}

\newtheorem{conj}[prop]{Conjecture}
\newtheorem{rem}[prop]{Remark}
\newtheorem{coro}[prop]{Corollary}

\title{Towards the Mirror Symmetry for Calabi-Yau Complete Intersections
in  Gorenstein Toric Fano Varieties}

\author{Lev A. Borisov \\
Department of Mathematics \\
University of Michigan \\
Ann Arbor, Michigan 48109-1003 \\
e-mail: Lev.Borisov@math.lsa.umich.edu}

\begin{document}

\date{}

\maketitle

\begin{abstract}
We propose  a  combinatorical duality
for lattice polyhedra  which conjecturally gives rise
to the pairs of mirror symmetric families of Calabi-Yau complete
intersections in toric Fano varieties with Gorenstein singularities.
Our construction is a generalization of the polar duality proposed
by Batyrev for the case of hypersurfaces.
\end{abstract}

\section{Introduction}

\noindent

Mirror Symmetry discovered by physicists for Calabi-Yau manifolds
still remains a surprizing puzzle for mathematicians. Some insight on this
phenomenon was received
from the investigation of
Mirror Symmetry for some examples of  Calabi-Yau varieties which
admit simple birational models embedded in toric varieties.
In this context, Calabi-Yau manifolds obtained by the resolution
of singularities of complete intersections in toric varieties are
the most general examples.

In the  paper of Batyrev and van Straten \cite{bat.strat},
there was proposed a method for conjectural construction
of mirror families for Calabi-Yau complete intersections
in  toric varieities. Unfortunately, their method fails to provide such a nice
duality as it is in the case of hypersurfaces \cite{bat.dual}. The purpose of
these
notes is to propose a generalized  duality which  conjecturally gives rise
to the mirror involution for complete intersections.

I am pleased to thank prof. Batyrev who has edited my original notes.

\section{Basic definitions and notations}

\noindent

Let $M$ and $N = {\rm Hom}(M, {\bf Z})$ be
dual free abelian groups of rank $d$,
$M_{\bf R}$ and $N_{\bf R}$ be their real scalar extensions and
\[ \langle \cdot , \cdot \rangle\;  :\;
  M_{\bf R} \times N_{\bf R} \rightarrow {\bf R} \]
be the canonical pairing. For any convex polyhedron $P$ in
$M_{\bf R}$ (or in $N_{\bf R})$, we denote its set of vertices by $P^0$.

\begin{dfn}
{\rm Let $P$ be a $d$-dimensional
convex  polyhedron in $M_{\bf R}$ such
that $P$ contains zero point $0 \in M_{\bf R}$ in its interior.
Then
\[ P^* = \{ y \in N_{\bf R} \mid \langle x, y \rangle \geq -1 \} \]
is called {\em polar}, or {\em dual polyhedron}.  }
\label{polar}
\end{dfn}

\begin{dfn}
{\rm A convex polyhedron $P$ in $M_{\bf R}$ is called a {\em lattice
polyhedron} if $P^0 \subset M \subset M_{\bf R}$. }
\end{dfn}

\begin{dfn}
{\rm (cf. \cite{bat.dual})  Let $\Delta$ be a $d$-dimensional
lattice polyhedron in $M_{\bf R}$ such
that $\Delta$ contains $0$ in its interior.
Then $\Delta$ is called {\em reflexive} if
$\Delta^*$ is also a lattice polyhedron.}
\end{dfn}

\begin{dfn}
{\rm Let $P$ be a $d$-dimensional
convex  polyhedron in $M_{\bf R}$ such
that $P$ contains zero point $0 \in M_{\bf R}$ in its interior.
We define the $d$-dimensional fan $\Sigma \lbrack P \rbrack$ as
the union of the zero-dimensional cone $\{ 0 \}$ together with
the set of all cones
\[ \sigma\lbrack \theta \rbrack = \{ 0 \} \cup
\{ x \in M_{\bf R} \mid
\lambda x \in \theta\; \mbox{for some $\lambda \in {\bf R}_{> 0}$}  \}  \]
supporting faces $\theta$ of $P$. }
\end{dfn}

Next four definitions of this section
play the main role in our construction.

\begin{dfn}
{\rm Let $\Delta \in M_{\bf R}$ be a reflexive polyhedron. Put $E =
\{ e_1, \ldots, e_n \} =
\Delta^0$. A representation  of $E = E_1 \cup \cdots \cup E_r$
as the union of disjoint subsets $E_1, \ldots , E_r$ is called
{\em nef-partition of} $E$ if there exist integral convex
$\Sigma \lbrack \Delta \rbrack $-piecewise linear functions
$\varphi_1, \ldots, \varphi_r $ on $M_{\bf R}$ such that
$\varphi_i(e_j) = 1$ if $ e_j \in E_i$, and $\varphi_i(e_j) = 0$ otherwise.}
\end{dfn}

\begin{rem}
{\rm The term {\em nef-partition} is motivated by the fact that
such a partition induces a representation of the anticanonical
divisor $-K$ on the Gorenstein toric Fano variety ${\bf P}_{\Delta^*}$
as the sum of $r$ Cartier divisors which are nef. }
\end{rem}

\begin{dfn}
{\rm Let $E = E_1 \cup \cdots \cup E_r$ be a nef-partition. Define $r$
convex polyhedra $\Delta_1, \ldots, \Delta_r \subset M_{\bf R}$ as
\[ \Delta_i = {\rm Conv}(\{0\} \cup   E_i ), \;
i =1, \ldots, r. \]}
\label{delta}
\end{dfn}

\begin{rem}
{\rm From Definition \ref{delta}  we immediately obtain that
$\Delta_i \cap  \Delta_j = \{0 \}$ if $i \neq j$ and
$\Delta =
{\rm Conv}( \Delta_1 \cup \cdots \cup \Delta_r )$.  }
\end{rem}

\begin{dfn}
{\rm Let $E = E_1 \cup \cdots \cup E_r$ be a nef-partition. Define $r$
convex polyhedra $\nabla_1, \ldots, \nabla_r \subset N_{\bf R}$ as
\[ \nabla_i = \{ y \in N_{\bf R} \mid \langle x, y \rangle \geq - \varphi_i(x)
\},
\; i =1, \ldots, r. \] }
\label{nabla}
\end{dfn}

\begin{rem}
{\rm It is obvious that $\{ 0 \} \in \nabla_1 \cap \cdots \cap \nabla_r$.
By Definition \ref{polar},  one has
\[ \Delta^* = \{ y \in N_{\bf R} \mid \langle x, y \rangle \geq - \varphi(x)
\},  \]
where $\varphi  = \varphi_1 + \cdots \varphi_r$. Therefore
$\nabla_1 \cup \cdots \cup \nabla_r \subset \Delta^*$.
Notice that $\nabla_1, \ldots, \nabla_r$ are also lattice polyhedra.
This fact follows from the following standard statement.}
\label{prop.nabla}
\end{rem}

\begin{prop}
Let $\Sigma$ be any complete fan of cones in $M_{\bf R}$, $\varphi_0$  a
convex $\Sigma$-piecewise linear function on $M_{\bf R}$. Then
\[ Q_0 =  \{ y \in N_{\bf R} \mid \langle x, y \rangle \geq -
\varphi_0(x) \} \]
is a convex polyhedron whose vertices are restrictions of $\varphi_0$ on
cones  of maximal dimension of $\Sigma$.
\label{bas.lem}
\end{prop}

\begin{coro}
The convex functions $\varphi_1, \ldots, \varphi_r$ have form
\[ \varphi_i(x) = - \min_{y \in \nabla_i} \langle x , y \rangle. \]
In particular, we have
\[ - \min_{  x \in \Delta_j^0,\, y \in \nabla_i^0}
 \langle x, y \rangle = \delta_{j\,i} \]
and
\[ \langle \Delta_j, \nabla_i \rangle \geq - \delta_{j\,i}. \]
\label{relations}
\end{coro}

\begin{dfn}
{\rm Define the lattice polyhedron $\nabla \in N_{\bf R}$ as
\[ \nabla = {\rm Conv}(\nabla_1 \cup \cdots \cup \nabla_r). \]}
\end{dfn}

\begin{rem}
{\rm Remark \ref{prop.nabla} shows that $\nabla \subset \Delta^*$.}
\label{inclusion}
\end{rem}

\section{The combinatorical duality}

\begin{prop}
$\Delta^* = \nabla_1 + \cdots +  \nabla_r$.
\end{prop}

{\em Proof.} The statement follows from the equality $\sum_i \varphi_i =
\varphi$,
from Remark \ref{prop.nabla} and Proposition \ref{bas.lem}.  \hfill
$\Box$

\begin{prop}
$\nabla^* = \Delta_1 + \cdots + \Delta_r$.
\end{prop}

{\em Proof. }  Let $x = x_1 + \cdots + x_r$ be a point of
$\Delta_1 + \cdots + \Delta_r$ $(x_i \in \Delta_i)$, and $y = \lambda_1 y_1 +
\cdots \lambda_r y_r$, $(\lambda_1 + \cdots +
\lambda_r = 1,\;  \lambda_i \geq 0, \; y_i \in \nabla_i)$
be a point in $\nabla$. By
\ref{relations},
\[ \langle x, y \rangle \geq \sum_{i=1}^r \lambda_i \langle x_i, y_i \rangle
\geq - \sum_{i=1}^r  \lambda_i = -1. \]
Hence  $\Delta_1 + \cdots + \Delta_r \subset \nabla^*$.

Let $y \in (\Delta_1 + \cdots + \Delta_r)^*$. Put
\[ \lambda_i = - \min_{ x \in \Delta_i} \langle x, y \rangle. \]
Since $0 \in \Delta_i$, all $\lambda_i$ are nonnegative.
Since $ \langle \sum_i \Delta_i , y \rangle \geq -1$, we have
$\sum_i \lambda_i \leq 1$. Consider the convex function
$\varphi_y = \sum_i \lambda_i \varphi_i$. For all $x \in M_{\bf R}$, we have
\[ - \varphi_y(x) =
\sum_{i=1}^r  \lambda_i \varphi_i(x) \leq \langle x, y \rangle . \]
By Proposition \ref{bas.lem}, $y$ is contained in the convex hull of
all points in $N_{\bf R}$ which are equal to restrictions
of $\varphi_y$ on cones of maximal dimension of $\Sigma
\lbrack \Delta \rbrack$.
By definition of $\varphi_y$, any such a point is a sum
$\sum_{i} \lambda_i p_i$ where $p_i \in \nabla_i$. Hence
$y \in \nabla$. Thus we have proved that
$(\Delta_1 + \cdots + \Delta_r)^* \subset \nabla$. \hfill $\Box$

Since $\nabla$ and $\Delta_1 + \cdots + \Delta_r$ are lattice polyhedra,
we obtain:

\begin{coro}
The polyhedron $\nabla$ is reflexive.
\end{coro}

\begin{prop}
Let $E'= \{e_1', \ldots , e_k' \} = \nabla^0$, $E_i' = \nabla^0_i$
$( i =1, \ldots r)$. Then subsets
$E_1',  \ldots E_r' \subset E'$ give rise to a  nef-partition  of $E'$.
\end{prop}

{\em Proof. } First, we prove that $\nabla_i \cap \nabla_j = \{ 0 \}$ for
$i \neq j$. Assume that $e_p' \in \nabla_i \cap \nabla_j$.
Using \ref{relations},  we obtain that $e_p'$ has non-negative
values at all vertices
$e_1, \ldots, e_n$ of $\Delta$. On the other hand, $e_p'$ has zero
value at the interior point $0 \in \Delta$. Hence  $e_p'$ must be
zero. This means that $E_i' \cap E_j' = \emptyset$ for $i \neq j$.

Let $e_p'$ be a vertex of $\nabla_i$. We prove that $e_p'$ is also a
vertex of $\nabla$. By \ref{relations},  there exists
a vertex $e_s \in \Delta_j^0$
such that $\langle e_s, e_p'\rangle = -1$. Moreover,
\[ -1 =\min_{y \in \nabla} \langle e_s, y \rangle =
\min_{y \in \nabla_i^0} \langle e_s, y \rangle. \]
So $e_p'$ is also a vertex of $\nabla$.

Define the functions
\[ \psi_i \; : \; N_{\bf R} \rightarrow {\bf R},\; i =1, \ldots, r; \]
\[ \psi_i(y) =  - \min_{x \in \Delta_i}  \langle x, y \rangle. \]
Obviously, $\psi_1, \ldots, \psi_r$ are convex. By \ref{relations},
$\psi_i(e_p') = 1$ if $e_p' \in \nabla_i$, and $\psi_i(e_p') = 0$ otherwise.
We prove that  restrictions of $\psi_i$ on cones of
$\Sigma \lbrack \nabla \rbrack$ are linear. It is sufficient to consider
restrictions of $\psi_i$ on cones $\sigma \lbrack \theta \rbrack$
of maximal dimension where $\theta = \nabla \cap
\{y \mid \langle v, y \rangle = -1\}$ is a $(d-1)$-dimensional face
of $\nabla$ corresponding to a vertex $v \in \nabla^* =
\Delta_1 + \cdots + \Delta_r$. Let
$v = v_1 + \cdots + v_i +\cdots + v_r$,  where  $v_i$ denotes
a vertex of $\Delta_i$.  If we take another vertex $v_i' \neq v_i$ of
$\Delta_i$, then the sum $v = v_1 + \cdots + v_i' +\cdots + v_r$
represents another vertex of $\nabla^*$. Clearly,
$\langle v, y \rangle \leq \langle v', y \rangle$ for any
$y \in \sigma\lbrack \theta \rbrack$, i.e.,
$\langle v_i, y \rangle \leq \langle v_i', y \rangle$. Hence
the restriction of $\psi_i$ on $\sigma\lbrack \theta \rbrack$
is $- \langle v_i, y \rangle$. \hfill $\Box$

\begin{coro}
\[ \Delta_i = \{ x \in m_{\bf R} \mid \langle x, y \rangle \geq - \psi_i(y) \},
\; i =1, \ldots, r. \]
\end{coro}
\medskip

Thus we have proved that the set of reflexive polyhedra with
nef-partitions has a natural involution
\[ \imath\; : \; (\Delta; E_1, \ldots, E_r) \rightarrow
(\nabla; E_1', \ldots, E_r' ). \]
On the other hand, every nef-partition of a reflexive polyhedron
$\Delta$ defines $r$ base point free linear systems of numerically effective
Cartier divisors $\mid D_1 \mid , \ldots,
\mid D_r \mid$ such that the sum $D_1 + \ldots + D_r$ is
the anticanonical divisor on the Gorenstein toric Fano variety
${\bf P}_{\Delta^*}$.

\begin{conj}
The duality between nef-partitions of reflexive polyhedra $\Delta$ and
$\nabla$ gives  rise to pairs of mirror symmetric
families of Calabi-Yau complete intersections
in Gorenstein toric Fano  varieties ${\bf P}_{\Delta^*}$ and
${\bf P}_{\nabla^*}$.
\end{conj}

\end{document}